\documentclass[12pt]{amsart}

\usepackage{amsthm}

\usepackage{mathrsfs,mathtools,latexsym,eucal}
\usepackage{amscd,amsfonts,amssymb,amsmath,amsthm}
\usepackage{mathrsfs}
\usepackage{graphicx,graphics,color}
\usepackage[utf8]{inputenc}
\usepackage[active]{srcltx}
\usepackage[english]{babel}
\usepackage[pagewise,displaymath,mathlines]{lineno}
\usepackage{listings}
\usepackage{color}
\usepackage{subfigure} 
\usepackage{url}
\usepackage{float}
\usepackage{enumerate}

\usepackage[active]{srcltx}

\usepackage{pgf,tikz}
\usetikzlibrary{arrows}

\definecolor{dkgreen}{rgb}{0,0.6,0}
\definecolor{gray}{rgb}{0.5,0.5,0.5}
\definecolor{mauve}{rgb}{0.58,0,0.82}

\newtheorem*{theorem*}{Theorem}
\newtheorem{theorem}{Theorem}[section]

\newtheorem{proposition}[theorem]{Proposition}
\newtheorem{corollary}[theorem]{Corollary}

\newtheorem{remark}[theorem]{Remark}

\newtheorem{example}[theorem]{Example}

\renewcommand{\leq}{\leqslant}
\renewcommand{\geq}{\geqslant}

\numberwithin{equation}{section}

\newcommand{\card}{\mathrm{card}}

\newcommand{\disp}{\displaystyle}

\begin{document}

\title{Completeness of a normed space via strong $p$-Cesàro summability}


\author{Fernando León-Saavedra}
\address{Departamento de Matemáticas, Universidad de Cádiz, Facultad de Ciencias Sociales y de la Comunicación, Avenida de la Universidad s/n, 11403-Jerez de la Frontera (Cádiz)}
\email{{\tt fernando.leon@uca.es}}

\author{Soledad Moreno-Pulido}
\address{Department of Mathematics, College of Engineering, University of Cadiz, Puerto Real 11510, Spain (EU)}
\email{{\tt soledad.moreno@uca.es}}

\author{Antonio Sala-Pérez}
\address{Department of Mathematics, College of Engineering, University of Cadiz, Puerto Real 11510, Spain (EU)}
\email{{\tt antonio.sala@uca.es}}

\begin{abstract}
In this paper we will characterize the completeness and barrelledness of a normed space through the strong $p$-Cesáro convergence of series. A new characterization of weakly unconditionally Cauchy series and unconditionally convergent series through the strong $p$-Cesàro summability is obtained. 
\end{abstract}


\maketitle \thispagestyle{empty}






\section{Introduction}
Let $X$ be a normed space and $0<p<\infty$, a sequence $(x_k)$ is said to be strongly p-Cesàro convergent to $L\in X$ if
$$
\lim_{n\to \infty} \frac{1}{n}\sum_{k=1}^n \|x_k-L\|^p=0.
$$
The strong  1-Cesàro convergence for real numbers was introduced by Hardy-Littlewood \cite{hardy} and Fekete \cite{fekete} in connection with the convergence of Fourier series (see \cite{zeller-beekmann}, for historical notes, and the most recent monograph \cite{boos}). 

Some years later, in 1935, Professor A. Zygmund (see \cite{zygmund} for one of the reprints) introduced the idea of statistical convergence in a independently way. A sequence $(x_n)$ is statistically convergent  to $L$ if for any $\varepsilon>0$ the subset $\{k\,:\, \|x_k-L\|<\varepsilon\}$ has density $1$ on the natural numbers. 

Both concepts  were developed independently and surprisingly enough, both are related thanks to a result by J. Connor (\cite{Connor88}).
Since then, in this circle of ideas, a significant number of deep and beautiful results have been obtained by Connor, Fridy, Mursaleen...and many others (see \cite{b3,b4,b5,b2,b1,b7,b6,b8})

There are also results that obtain characterizations of properties of Banach spaces through convergence types. For instance, Kolk  \cite{kolk} was one of the pionnering contributors. Connor, Ganichev and Kadets  \cite{connorganichevkadets} obtained important results that
relate the statistical convergence to classical properties of Banach spaces.

The aim of this paper is to obtain properties of a Banach space studying properties of strong $p$-Cesàro convergence of a series. 
Let $X$ be a normed space, and set $\sum x_i$  a series in $X$.  In \cite{Aizpuru2000} the spaces of convergence $S(\sum x_i)$ associated to the series $\sum x_i$ are introduced. $S(\sum x_i)$ is defined as the sequences $(a_j)\in \ell_\infty$ such that $\sum a_ix_i$ converges.
The space $X$ is complete if and only if for every weakly unconditionally Cauchy series $\sum x_i$, the space $S(\sum x_i)$ is complete. Moreover, the space $X$ is barrelled if and only if each series $\sum_i x_i^\star$ in $X^\star$ if the corresponding space of weak-$\star$ convergence associated to $\sum x_i^\star$ is the entire space $\ell_\infty$, that is, $S_{w^\star}(\sum x_i^\star)=\ell_\infty$.

In this paper we explore this structure for the strong-$p$ Cesàro convergence. At first glance, it seems that in order to show that a sequence is strongly $p-$Cesàro convergent it is necessary to know the value of its limit previously. However, thanks to the results by Connor \cite{Connor88} and Fridy \cite{Fridy}, we can avoid this difficulty. Section 2 is an expository session where we will show examples and preliminary aspects related to the strong  $p$-Cesàro convergence.
Section 3 deals with space of Strog-Cesàro convergence. It is shown that a series in a Banach space is weakly unconditionally Cauchy if and only if its space of strong $p$-Cesàro convergence is complete. Moreover, if this equivalence is true for each series in a normed space $X$, then the space $X$ must be complete.
In Section 4 and 5 we will begin by defining reasonably, the strong  $p$-Cesàro convergence for the weak and the weak-$\star$ topology in a Banach space $X$ and its dual $X^\star$ respectively. After this, we will show analogous results for the strong $p$-Cesàro convergence in these topologies. We also prove a characterization of barrelledness which is similar to the aforementioned one, but replacing weak-$\star$ convergence by our concept of strong $p$-Cesàro convergence for the weak-$\star$ topology.

\section{Some preliminary results}

We begin this section by recalling some preliminaries we will need throughout this work. If $A\subset\mathbb{N}$, the density of $A$ is denoted by 
    $
    d(A)=\lim_n\frac{1}{n}\card(\{k\leq n: k\in A\}),
    $
    whenever this limit exists.

Let $X$ be a normed space and $x=(x_k)_k$ a sequence in $X$. The sequence $x$ is said to be \emph{statistically convergent} if there is $L\in X$ such that for every $\varepsilon>0$, 
$
	d\bigl(\{k: \|x_k-L\|\geq\varepsilon\}\bigr)=0
$
or equivalently
$
	d\bigl(\{k: \|x_k-L\|<\varepsilon\}\bigr)=1
$
and we will write $(x_k)\overset{st}{\to}L$ and $L=st-\lim_n x_n$.
The sequence $x$ is said to be \emph{statistically Cauchy} if for each $\varepsilon>0$ and $n\in\mathbb{N}$, there exists $p\geq n$ such that 
$
d\bigl(\{k: \|x_k-x_p\|\geq\varepsilon\}\bigr)=0
$ 
or equivalently
$
d\bigl(\{k: \|x_k-x_p\|<\varepsilon\}\bigr)=1.
$ 

Fridy \cite[Theorem 1]{Fridy} proved that in a Banach space, a sequence is statistically convergent if and only if it is statistically Cauchy.
	
Let us consider now $0<p<+\infty$. The sequence $x$ is said to be \emph{strongly $p-$Cesàro or} $\textrm{w}_p$ \emph{summable} if there is $L\in X$ such that
$
	\displaystyle\lim_n\frac{1}{n}\sum_{k=1}^n\|x_k-L\|^p=0, 
$
in which case we say that $x$ is strongly $p-$Cesàro summable to $L$, and we will write $(x_k)\overset{\textrm{w}_p}{\to}L$ and $L=\textrm{w}_p-\lim_n x_n$.

Although the convergent sequences are $\textrm{w}_p$ convergent, it is easy to see that this kind of convergence is weaker than the usual, as we will show in the next example:

\begin{example}

There exist unbounded sequences that are strong $p-$Cesàro summable.

Let us consider the sequence $n_j=j^3$ for every $j\in\mathbb{N}$ and define
$$
	x_k=\left\{
    	\begin{array}{ll}
    		0, & k\neq r^3\mbox{ for all }r.\\
            r, & k=r^3\mbox{ for some }r.
        \end{array}
    \right.
$$
The sequence $(x_k)$ is unbounded. For every $n\in\mathbb{N}$, let $n_r=\max\{n_j:n_j\leq n\}$. Then:
	\begin{align*}
		\frac{1}{n}\sum_{k=1}^n |x_k| = \frac{1}{n}\sum_{k=1}^n x_k&\leq \frac{1}{n_r}\sum_{k=1}^{n_r}x_k=\frac{1+\dots+r}{r^3}\underset{r\to\infty}{\to} 0.
	\end{align*}
\end{example}

However, a sequence that is $\textrm{w}_p$ summable which is unbounded cannot diverge randomly, as the following proposition shows:

\begin{proposition}\label{propo}
Let $0<p<\infty$ and $(x_k)_k$ be a sequence in a normed space $X$ such that for some increasing subsequence $(n_j)\subset \mathbb{N}$,  $\displaystyle\lim_j\frac{1}{n_j}\sum_{k=1}^{n_j}\|x_k\|^p=+\infty$. Then, $(x_k)_k$ is not strongly $p-$Cesàro summable to any $L\in X$.
\end{proposition}
{\it Proof. } Suppose on the contrary that there exists $L\in X$ such that $$\lim_n\frac{1}{n}\sum_{k=1}^n\|x_k-L\|^p=0.$$ Then:
    \begin{align*}
    	\frac{1}{n_j}\sum_{k=1}^{n_j} \|x_k\|^p & \leq \frac{1}{n_j}\sum_{k=1}^{n_j} (\|x_k-L\|+\|L\|)^p \\
      & \leq \frac{1}{n_j}\sum_{k=1}^{n_j} \|x_k-L\|^p+\|L\|^p 
    \end{align*}
    which converges to $\|L\|^p$ as $n_j\to \infty$, which is a contradiction because the first part of the inequality diverges by hypothesis.$\hfill\Box$

Connor \cite[Theorem 2.1]{Connor88} discovered that the real bounded sequences $\textrm{w}_p$ convergent are exactly the  statistically convergent sequences. This fact also holds for normed spaces and we include the proof for the sake of completeness.

\begin{proposition}[Connor \cite{Connor88}]\label{Connor_normed}
	Set $0<p<\infty$ and let $X$ be a normed space. If a sequence is strongly $p-$Cesàro summable to $L$, then it is statistically convergent to $L$. Additionally, if the sequence is bounded, the converse is also true.
\end{proposition}
{\it Proof. } Let us consider $(x_k)_k$ a sequence which is strongly $p-$Cesàro summable to $L\in X$ and $\varepsilon>0$. For any $n\in\mathbb{N}$,
    \begin{align*}
    	\sum_{k=1}^n\|x_k-L\|^p\geq \sum_{\substack{k=1 \\ \|x_k-L\|\geq\varepsilon}}^n\|x_k-L\|^p\geq \sum_{\substack{k=1 \\ \|x_k-L\|\geq\varepsilon}}^n \varepsilon^p=\card\bigl(\{k\leq n: \|x_k-L\|^p\geq\varepsilon\}\bigr)\;\varepsilon^p.
    \end{align*}
Since $(x_k)_k$ is strongly $p-$Cesàro summable to $L$, we have that $\displaystyle\lim_n\displaystyle\frac{1}{n}\displaystyle\sum_{k=1}^n||x_k-L||^p=0$, so for every $\varepsilon>0$, $\displaystyle\lim_n\frac{1}{n}\card\{k\leq n: ||x_k-L||\geq\varepsilon\}=0$  which shows that $(x_k)_k$ is statistically convergent to $L$.

Suppose now that $x=(x_k)_k$ is a bounded sequence which is statistically convergent to $L\in X$ and set $K=\|x\|_{\infty}+\|L\|$, where $\|x\|_{\infty}=\displaystyle\sup_k \|x_k\|$. Given $\varepsilon>0$, there exists $N_{\varepsilon}\in\mathbb{N}$ such that 
\begin{align*}
	\frac{1}{n}\card\left(\Bigl\{k\leq n: \|x_k-L\|\geq \left(\frac{\varepsilon}{2}\right)^{1/p}\Bigr\}\right)<\frac{\varepsilon}{2 K^p},
\end{align*}
for every $n\geq N_{\varepsilon}$. Set $L_n=\{k\leq n:\|x_k-L\|\geq \left(\frac{\varepsilon}{2}\right)^{1/p}\}$. For every $n\geq N_{\varepsilon}$, we have:
\begin{align*}
	\frac{1}{n}\sum_{k=1}^n\|x_k-L\|^p & =\frac{1}{n}\left[\sum_{k\in L_n} \|x_k-L\|^p+\sum_{\substack{k\leq n \\ k\notin L_n }}\|x_k-L\|^p\right] \\
                                       & <\frac{1}{n} K^p \card(L_n)+\frac{1}{n} n\frac{\varepsilon}{2} \\
                                       & < \frac{1}{n} K^p n\frac{\varepsilon}{2 K^p} + \frac{\varepsilon}{2}=\varepsilon.
\end{align*}
Thus, $(x_k)_k$ is strongly $p-$Cesàro summable to $L$.$\hfill\Box$

Next, we show that for the converse, boundedness is necessary.

\begin{example}

There exist unbounded  statistically convergent sequences which are not strongly $p-$Cesàro summable. Indeed, set $n_j=j^2$  and let us define
$$
	x_k=\left\{
    	\begin{array}{ll}
    		0, & k\neq j^2\mbox{ for all }j.\\
            j^{2/p}, & k=j^2\mbox{ for some }j.
        \end{array}
    \right.
$$
The sequence $(x_k)_k$ is unbounded. Take $\varepsilon>0$, it is easy to see that $\card\bigl(\{k\leq n: |x_k|\geq\varepsilon\}\bigr)=0$, so $(x_k)_k$ is statistically convergent to zero.Let us Observe that:
\begin{align*}
		\frac{1}{n_j}\sum_{k=1}^{n_j} |x_k|^p = \frac{1}{n_j}\sum_{k=1}^{n_j} x_k^p=\frac{1+2^2+3^2+\dots+j^2}{j^2},
	\end{align*}
which diverges as $j\to \infty$. Hence, by applying Proposition \ref{propo}, we deduce that $(x_k)_k$ is not strongly $p-$Cesàro summable.
\end{example}

Let us recall that a sequence $x=(x_k)_k$ in a normed space $X$ is said to be \emph{Cesàro convergent} if there is $L\in X$ such that
 $\displaystyle\lim_n \left\|\frac{1}{n}\sum_{k=1}^n x_k-L\right\|=0$. The $\textrm{w}_p$ summability is related to the Cesàro convergence in a natural way:

\begin{proposition}
\label{propoantonio}
Let $X$ be a normed space and $(x_k)_k$ a sequence in $X$. If $p\geq1$ and $(x_k)_k$ is strongly $p-$summable to $L$, then $(x_k)_k$ is Cesàro convergent to $L$. 

\end{proposition}
{\it Proof. } Let us observe that
\begin{align*}
	0 &\leq \left\|\frac{1}{n}\sum_{k=1}^n x_k -L\right\|=\frac{1}{n}\left\|\sum_{k=1}^n x_k-nL\right\|=\frac{1}{n}\left\|\sum_{k=1}^n(x_k-L)\right\| \\
     &\leq \frac{1}{n}\sum_{k=1}^n \|x_k-L\|\leq \frac{1}{n}\sum_{k=1}^n \|x_k-L\|^p\underset{n\to\infty}{\to}  0
\end{align*}
$\hfill\Box$

\begin{remark}
Let us observe that the condition $p\geq 1$ is sharp. Indeed, the sequence

$$
	x_k=\left\{
    	\begin{array}{ll}
    		0, & k\neq r^3\mbox{ for all }r.\\
            r^2, & k=r^3\mbox{ for some }r.
        \end{array}
    \right.
$$
is $\frac{1}{2}-$ Cesàro summable to zero, and the Cesàro means do not converge to zero.
\end{remark}

The converse of Proposition \ref{propoantonio} is clearly not true, as we show in the next example.

\begin{example} There exist  Cesàro convergent sequences which are not $p-$Cesàro summable. Let us define
$$
	x_k=\left\{\begin{array}{rl}
    			 -1 & \mbox{ if } k \mbox{ is odd}, \\
                  0 & \mbox{ if } k \mbox{ is even}.
    		   \end{array}
    	\right.
$$

$(x_k)_k$ is not strong $p-$Cesàro summable to any $L\in\mathbb{R}$ because it is not statistically convergent to any $L$. However, observe that:
$$
	\left|\frac{1}{n}\sum_{k=1}^n x_k\right|=\frac{n/2}{n}\to \frac{1}{2},
$$
so $(x_n)_n$ is Cesàro convergent to $\frac{1}{2}$.
\end{example}
Finally, for future references, by the Stolz-Cesàro Theorem, we get:
\begin{proposition}\label{propoStolz}
	Observe that, if $(x_k)_k$ is a sequence in $\mathbb{R}$ and $\displaystyle\sum_{k=1}^{\infty}x_k=L$ where $L\in\mathbb{R}\cup\{\pm\infty\}$, then $\displaystyle\lim_n \frac{1}{n}\sum_{k=1}^n S_k=L$.
\end{proposition}

\section{The  strong $p-$Cesàro summability space}

Let $\sum_i x_i$ be a series in a real Banach space $X$, set $0<p<+\infty$ and let us define 
$$
	S_{\textrm{w}_p}\left(\sum_i x_i\right)=\left\{(a_i)_i\in\ell_{\infty}: \sum_{i} a_i x_i\mbox{ is }\textrm{w}_p\mbox{ summable}\right\}
$$
endowed with the supremum norm. This space will be called the space of $\textrm{w}_p$ summability associated to the series $\sum_i x_i$. The following theorem characterizes the completeness of the space $S_{\textrm{w}_p}\bigl(\sum_i x_i\bigr)$.

\begin{theorem}\label{Swpcompleto}
	Let $X$ be a real Banach space and $0<p<+\infty$. The following conditions are equivalent:
    \medskip
    \begin{enumerate}[(1)]
    	\item $\sum_i x_i$ is a weakly unconditionally Cauchy series (wuc).  
    	\item \emph{$S_{\textrm{w}_p}(\sum_i x_i)$} is a complete space. 
        \item \emph{$c_0\subset S_{\textrm{w}_p}(\sum_i x_i).$}
    \end{enumerate}
\end{theorem}

{\it Proof. } Let us show that (1)$\Rightarrow$(2). Since $\sum x_i$ is wuc, the following supremum is finite:
$$
	H=\sup\left\{\left\|\sum_{i=1}^n a_i x_i\right\|:|a_i|\leq 1,1\leq i\leq n,n\in\mathbb{N}\right\}<+\infty.
$$

Let $(a^m)_m\subset S_{\textrm{w}_p}(\sum_i x_i)$ such that $\displaystyle\lim_m\|a^m-a^0\|_{\infty}=0$, with $a^0\in\ell_{\infty}$. We will prove that $a^0\in S_{\textrm{w}_p}(\sum_i x_i)$. Let us suppose without any loss of generality that $\|a^0\|_{\infty}\leq 1$. Then, the partial sums $S_k^0=\sum_{i=1}^k a_i^0 x_i$ satisfy $\|S_k^0\|\leq H$ for every $k\in\mathbb{N}$, that is, the sequence $(S_k^0)$ is bounded. Then, $a^0\in S_{\textrm{w}_p}(\sum_i x_i)$ if and only if $(S_k^0)$ is $\textrm{w}_p$ summable to some $L\in X$. Since $(S_k^0)$ is bounded, according to Connor's Theorem \cite[Theorem 2.1]{Connor88} (Proposition
\ref{Connor_normed}), $(S_k^0)$ is $\textrm{w}_p$ summable if and only if $(S_k^0)$ is statistically convergent to some $L\in X$. According to \cite[Theorem 1]{Fridy}, $(S_k^0)$ is statistically convergent to $L\in X$ if and only if $(S_k^0)$ is a statistically Cauchy sequence.

Set $\varepsilon>0$ and $n\in\mathbb{N}$. Then, we obtain statement (2) if we show that there exists $p_0\geq n$ such that 
$$
d\bigl(\{k: ||S_k^0-S^0_{p_0}||<\varepsilon\}\bigr)=1.
$$ 

Given $\varepsilon>0$, since $a^m\to a^0$ in $\ell_{\infty}$, there exists $m_0>n$ such that $\|a^m-a^0\|_\infty<\displaystyle\frac{\varepsilon}{4H}$ for all $m>m_0$, and since $S_k^{m_0}$ is statistically Cauchy, there exists $p_0\geq n$ such that the density $
\displaystyle d\left(\left\{k: ||S_k^{m_0}-S^{m_0}_{p_0}||<\frac{\varepsilon}{2}\right\}\right)=1.$ Fix $k$ such that 
\begin{align}\label{eq1}
\displaystyle\|S_k^{m_0}-S_{p_0}^{m_0}\|<\frac{\varepsilon}{2}.
\end{align}
We will show that $\|S_k^0-S_{p_0}^0\|<\varepsilon$, and this will prove that 
$$
	\left\{k:\|S_k^{m_0}-S_{p_0}^{m_0}\|<\frac{\varepsilon}{2}\right\}\subset\{k\,:\,\|S_k^0-S_{p_0}^0\|<\varepsilon\}.
$$
Since the first set has density 1, the second will also have density 1 and we will be done.

Let us observe first that for every $j\in\mathbb{N}$,
$$
	\left\|\sum_{i=1}^j\frac{4H}{\varepsilon}(a_i^p-a_i^{m_0})x_i\right\|\leq H,
$$
therefore
\begin{align}\label{eq2}
	\bigl\|S_j^0-S_j^{m_0}\bigr\|=\left\|\sum_{i=1}^j(a_i^0-a_i^{m_0})x_i\right\|\leq\frac{\varepsilon}{4}.
\end{align}
Then, by applying the triangular inequality,
\begin{align*}
	\bigl\|S_k^0-S_{p_0}^0\bigr\| & \leq \bigl\|S_k^0-S_k^{m_0}\bigr\|+\bigl\|S_k^{m_0}-S_{p_0}^{m_0}\bigr\|+\bigl\|S_{p_0}^0-S_{p_0}^{m_0}\bigr\|\\
    & <\frac{\varepsilon}{4}+\frac{\varepsilon}{2}+\frac{\varepsilon}{4}=\varepsilon.
\end{align*}
where the last inequality follows by applying \eqref{eq1} and \eqref{eq2}, which yields to the desired result.

Now, let us observe that if $S_{\textrm{w}_p}(\sum_i x_i)$ is a complete space, it contains the space of eventually zero sequences $c_{00}$ and therefore we get $(2)\Rightarrow(3)$.

Finally, let us show $(3)\Rightarrow (1)$. If the series $\sum x_i$ is not wuc, there exists $f\in X^{\ast}$ such that $\displaystyle\sum_{i=1}^{\infty}|f(x_i)|=+\infty$. Inductively, we will construct a sequence $(a_i)_i\in c_0$ such that $\sum_i a_i f(x_i)=+\infty$ and $a_i f(x_i)\geq 0$. If we denote by $\displaystyle S_k=\sum_{i=1}^k a_i f(x_i)$, then by applying Proposition \ref{propoStolz}, $\displaystyle\lim_n \frac{1}{n}\sum S_k=+\infty$. This implies that, by applying Proposition \ref{propo}, $(S_k)_k$ is not $\textrm{w}_p$ summable to any $L\in \mathbb{R}$, which is a contradiction with statement (3).

Since $\sum_{i=1}^{\infty}|f(x_i)|=+\infty$, there exists $m_1$ such that $\sum_{i=1}^{m_1}|f(x_i)|>2\cdot2$. We define $a_i=\frac{1}{2}$ if $f(x_i)\geq 0$ and $a_i=-\frac{1}{2}$ if $f(x_i)<0$ for $i\in\{1,2,\dots,m_1\}$. This implies that $\sum_{i=1}^{m_1}a_i f(x_i)>2$ and $a_i f(x_i)\geq 0$ if $i\in\{1,2,\dots,m_1\}$.

Let $m_2>m_1$ be such that $\sum_{i=m_1+1}^{m_2}|f(x_i)|>2^2\cdot 2^2$. We define $a_i=\frac{1}{2^2}$ if $f(x_i)\geq0$ and $a_i=-\frac{1}{2^2}$ if $f(x_i)<0$ for $i\in\{m_1+1,\dots,m_2\}$. Then, $\sum_{i=m_1+1}^{m_2}a_i f(x_i)>2^2$ and $a_i f(x_i)\geq 0$ if $i\in\{m_1+1,\dots,m_2\}$.

Inductively we obtain a sequence $(a_i)_i\in c_0$ with the above properties which lead us to a contradiction.$\hfill\Box$

\begin{remark}\label{nota3.1}

Let us observe that in the above proof, the completeness hypothesis is used in the implication $(1)\Rightarrow(2)$. Specifically, when we use Fridy's result (\cite[Theorem 1]{Fridy}). On the other hand, the implication $(2)\Rightarrow(3)$ that we will use in Theorem \ref{teor3.5} does not use the completeness of the space $X$.

\end{remark}

\begin{remark}

Let $\sum_i x_i$ be a series in a normed space $X$ and let
$$
	S\left(\sum_i x_i\right)=\left\{(a_i)_i\in \ell_{\infty}:\sum_i a_i x_i\mbox{ converges}\right\},
$$
endowed with the supremum norm. Clearly, $S\left(\sum_i x_i\right)$ is a subspace of $\ell_{\infty}$ and \emph{$S\left(\sum_i x_i\right)\subseteq S_{\textrm{w}_p}(\sum_i x_i)$}. If $X$ is a Banach space, then $\sum_i x_i$ is wuc if and only if $S(\sum_i x_i)$ is complete \cite{Aizpuru2000}. Theorem \ref{Swpcompleto} gives us a similar characterization by considering $\textrm{w}_p$ summability.
\end{remark}

\begin{corollary}

Let $X$ be a Banach space, $\sum_i x_i$ a series in $X$ and $p\geq 1$. The following properties are equivalent:

\begin{enumerate}[(1)]

	\item $\sum_i x_i$ is wuc.
    \item $S(\sum_i x_i)$ is a complete space.
    \item $c_0\subseteq S(\sum_i x_i)$.
    \item \emph{$S_{\textrm{w}_p}(\sum_i x_i)$} is a complete space.
    \item \emph{$c_0\subseteq S_{\textrm{w}_p}(\sum_i x_i)$}.
    \item \emph{$\sum|f(x_i)|$ is $\textrm{w}_p$} summable for every $f\in X^*$.

\end{enumerate}

\end{corollary}

{\it Proof. } The first three equivalence properties (1), (2) and (3) can be found in \cite{Aizpuru2000} and the rest of equivalences are  consequences of Theorem \ref{Swpcompleto}.$\hfill\Box$

Now let us show another main theorem.

\begin{theorem}\label{teor3.5}
Let $X$ be a normed space and $p\geq1$. Then $X$ is complete if and only if
\emph{$S_{\textrm{w}_p}(\sum_i x_i)$} is a complete space for every $\sum_i x_i$.
\end{theorem}
{\it Proof. } By Theorem \ref{Swpcompleto}, the condition is necessary. Now if $X$ is not complete, there exists $\sum x_i$ a series in $X$ such that $\|x_i\|\leq \frac{1}{i2^i }$ and $\sum x_i=x^{\ast\ast}\in X^{\ast\ast}\setminus X$. We will construct a wuc series $\sum_n y_n$ such that $S_{\textrm{w}_p}(\sum_n y_n)$ is not complete, a contradiction. 

Indeed, since $X^{\ast\ast}$ is a Banach space with the dual topology, $\disp\sup_{\|y^\ast\|\leq 1}|y^\ast(S_n)-x^{\ast\ast}(y^\ast)|\to0$, that is, $\disp\sum_{i=1}^\infty y^\ast(x_i)=x^{\ast\ast}(y^\ast)$, for all $\|y^\ast\|\leq 1$. By applying Proposition \ref{propoStolz}, we have:
\begin{align}\label{(A)}
\lim_{N\to\infty}\frac{1}{N}\sum_{k=1}^N y^\ast(S_k)=x^{\ast\ast}(y^\ast)
\end{align}
Set $y_n=n x_n$ and let us observe that $\sum y_n$ is a weakly unconditionally Cauchy series since $\|y_n\|<\frac{1}{2^n}$. We claim that the series $\sum_n \frac{1}{n}y_n$ is not $\textrm{w}_p$ summable in $X$.

On the contrary, let us suppose that $S_N=\sum_{n=1}^N \frac{1}{n}y_n$ is $\textrm{w}_p$ summable in $X$. That is, there exists $L\in X$ such that $\disp\lim_{N\to\infty}\frac{1}{N}\sum_{n=1}^N\|S_n-L\|^p=0$. In particular, for every $y^\ast\in X^\ast$ with $\|y^\ast\|\leq 1$ we have that $\disp\sup_{\|y^\ast\|\leq 1}\frac{1}{N}\sum_{k=1}^N|y^\ast(S_k-L)|^p\to0$. By applying Proposition \ref{propoantonio}, since $p\geq1$, we have that
\begin{align}\label{(B)}
\frac{1}{N}\sum_{k=1}^N y^\ast(S_k)=y^\ast(L), \mbox{ for every } \|y^\ast\|\leq 1.
\end{align}

From equations \ref{(A)} and \ref{(B)} and the uniqueness of the limit, we have that $x^{\ast\ast}(y^\ast)=y^\ast(L)$ for every $\|y^\ast\|\leq 1$, so we obtain $x^{\ast\ast}=L\in X$, which is a contradiction. This means that $S_N=\sum_{n=1}^N \frac{1}{n}y_n$ is not $\textrm{w}_p$ summable to any $L\in X$.

Finally, let us observe that since $\sum_n y_n$ is a weakly unconditionally Cauchy series and  $S_N=\sum_{n=1}^N\frac{1}{n}y_n$ is not $\textrm{w}_p$ summable, we have that $(\frac{1}{n})\notin \emph{$S_{\textrm{w}_p}(\sum_n y_n)$}$ and this means that $c_0\nsubseteq S_{\textrm{w}_p}(\sum_n y_n)$ which is a contradiction according to (3) in Theorem \ref{Swpcompleto} (see Remark \ref{nota3.1}) and the proof is completed. 
$\Box$

\begin{theorem}
	Let $\sum_i x_i$ be a series in a Banach space. The series $\sum_i x_i$ is wuc if and only if the operator \emph{$T:S_{\textrm{w}_p}(\sum_i x_i)\to X$} defined by \emph{$T((a_i)_i)=\textrm{w}_p-\lim_n S_n$  is continuous (where $S_n=\sum_{i=1}^n a_i x_i$}).
\end{theorem}

{\it Proof. } Suppose that $T$ is continuous and let us show that $\sum_i x_i$ is wuc. Since $c_{00}\subset S_{\textrm{w}_p}(\sum_i x_i)$, for every $(a_i)_i\in c_{00}$,
$
	\|T\bigl((a_i)_i\bigr)\|\leq \|T\| \|(a_i)_i\|_{\infty}.
$
Hence,
$$
	\sup_{n\in\mathbb{N}}\left\{\left\|\sum_{i=1}^n a_ix_i\right\|: |a_i|\leq 1\right\}\leq \|T\|,
$$
and this implies that the series $\sum x_i$ is wuc.

Let us suppose that $\sum_i x_i$ is wuc. Then, 
$$
	H=\sup_{n\in\mathbb{N}}\left\{\left\|\sum_{i=1}^n a_ix_i\right\|: |a_i|\leq 1\right\}<+\infty.
$$
Set $a=(a_i)_i\in S_{\textrm{w}_p}(\sum_i x_i)$ such that $\|a\|_{\infty}=1$. Then, $S_n=\sum_{i=1}^n a_i x_i$ is $\textrm{w}_p$ summable and since it is a bounded sequence, by applying Connor's Theorem \ref{Connor_normed}, it is statistically convergent to some $L$, and $L=st-\lim_n S_n = \textrm{w}_p-\lim_n S_n=T\bigl((a_i)_i\bigr)$. By applying Friddy's result \cite[Theorem 1]{Fridy}, there exists $A\subset \mathbb{N}$ of density 1 such that $\displaystyle\lim_{\substack{n \\ n\in A}}\|S_n-L\|=0$. For every $k\in A$, $\|S_k\|\leq H$, so
$$
	\|T\bigl((a_i)_i\bigr)\|=\|L\|=\lim_{\substack{k \\ k\in A}}\|S_k\|\leq \|H\|.
$$
which proves that $T$ is continuous and this completes the desired result.$\hfill\Box$

\section{The space of weak $\textrm{w}_p-$ summability}

In this section we study a similar structure with respect to the weak $\textrm{w}_p-$ summability.

	Let $X$ be a normed space. Set $0<p<+\infty$, a sequence $(x_k)_k$ is said to be \emph{weak} $\textrm{w}_p-$ \emph{summable} to $L\in X$ if for every $f\in X^{\ast}$, $f(x_k)\overset{\textrm{w}_p}{\to}f(L)$, that is,
$$
	\lim_n\frac{1}{n}\displaystyle\sum_{k=1}^n|f(x_k) - f(L)|^p=0, 
$$
and we will write $(x_k)\overset{w-\textrm{w}_p}{\to}L$ and $L=w-\textrm{w}_p-\lim_n x_n$.

Let $\sum_i x_i$ be a series in a Banach space $X$, $0<p<+\infty$. We now consider the space of $w-\textrm{w}_p$  summability given by:
$$
	S_{w-\textrm{w}_p}\left(\sum_i x_i\right)=\left\{(a_i)_i\in\ell_{\infty}: \sum_{i} a_i x_i\mbox{ is }w-\textrm{w}_p\mbox{ summable}\right\}
$$
endowed with the supremum norm. 

\begin{theorem}\label{Sw-wpcompleto}
	Let $0<p<+\infty$. The following conditions are equivalent:
    \medskip
    \begin{enumerate}[(1)]
    	\item $\sum_i x_i$ is a weakly unconditionally Cauchy series (wuc).  
    	\item \emph{$S_{w-\textrm{w}_p}(\sum_i x_i)$} is a complete space. 
        \item \emph{$c_0\subset S_{w-\textrm{w}_p}(\sum_i x_i).$}
    \end{enumerate}
\end{theorem}
{\it Proof. } Since $\sum_i x_i$ is wuc, 
    $$H=\displaystyle\sup\left\{\left\|\sum_{i=1}^n a_ix_i\right\|:|a_i|\leq 1,n\in\mathbb{N}\right\}<+\infty.$$ 
    Let $(a^m)_m\subset S_{w-\textrm{w}_p}(\sum_i x_i)$ be such that $\displaystyle\lim_m\|a^m-a^0\|_\infty=0$, with $a^0\in\ell_{\infty}$. We will prove that $a^0\in S_{w-\textrm{w}_p}(\sum_i x_i)$, and suppose without any loss of generality that $\|a^0\|_\infty\leq 1$. The sequence $S_k^0=\sum_{i=1}^k a_i^0x_i$ is bounded, and for every $f\in X^{\ast}$, we have that $f(S_k^0)=\sum_{i=1}^k a_i f(x_i)$ is a bounded sequence. We will show that $f(S_k^0)$ is $\textrm{w}_p$ summable. By applying again Connor's Theorem \ref{Connor_normed} and Friddy's result, it is sufficient to prove that $\bigl(f(S_k^0)\bigr)$ is statistically convergent or equivalently, $\bigl(f(S_k^0)\bigr)$ is statistically Cauchy. 
    
Given $\varepsilon>0$, we will show that for every $n\in\mathbb{N}$ there exists $p_0\geq n$ such that
$$
	d(\{ k: |f(S_k^0)-f(S_{p_0}^0)|\leq\varepsilon\})=1.
$$
Since $a^m\to a^0$ in $\ell_\infty$, there exists $m_0$ such that $\displaystyle\|a^m-a^0\|_\infty\leq \frac{\varepsilon}{4H\|f\|}$ for every $m\geq m_0$. And since $\bigl(f(S_k^{m_0})\bigr)$ is statistically Cauchy, for every $n\in\mathbb{N}$, there exists $p_0\geq n$ such that the set $\displaystyle\left\{k: |f(S_k^{m_0})-f(S_{p_0}^{m_0})|\leq\frac{\varepsilon}{2}\right\}$ has density 1. Let us consider $k\leq n$ such that 

\begin{align}\label{eq3}
	|f(S_k^{m_0})-f(S_{p_0}^{m_0})|<\frac{\varepsilon}{2}.
\end{align}

Let us observe that, for every $j$, $\displaystyle\left\|\sum_{i=1}^j\frac{\varepsilon}{4H\|f\|}(a_i^0-a_i^{m_0})x_i\right\|\leq H$, so we deduce that 

\begin{align}\label{eq4}
\displaystyle\|S_j^0-S_j^{m_0}\|=\left\|\sum_{i=1}^j(a_i^0-a_i^{m_0})x_i\right\|\leq \frac{\varepsilon}{4\|f\|}.
\end{align}

Then, using \eqref{eq3} and \eqref{eq4} and the triangular inequality,

\begin{align*}
	|f(S_k^0)-f(S_{p_0}^0)| & \leq |f(S_k^0 - S_k^{m_0})| + |f(S_{p_0}^0-S_{p_0}^{m_0})|+|f(S_k^{m_0}-S_{p_0}^{m_0})|\\
                            & \leq \|f\|\frac{\varepsilon}{4\|f\|} + \|f\|\frac{\varepsilon}{4\|f\|} + \frac{\varepsilon}{2}=\varepsilon,
\end{align*}
which implies that 
$$\displaystyle\left\{k: |f(S_k^{m_0})-f(S_{p_0}^{m_0})|\leq\frac{\varepsilon}{2}\right\}\subseteq \left\{k: |f(S_k^{0})-f(S_{p_0}^{0})|\leq\frac{\varepsilon}{2}\right\}$$ and since the first set has density 1, the second has also density 1 and we are done.

Finally, let us observe that implication $(2)\Rightarrow(3)$ is obvious and $(3)\Rightarrow (1)$ follows by a similar argument like in Theorem \ref{Swpcompleto}, and this finishes the proof.$\hfill\Box$

\section{The weak$^{\ast}$ $\textrm{w}_p-$ summability space}

We begin this section by defining a reasonable concept for weak$^{\ast}$ $\textrm{w}_p-$ summability. This convergence provides a different result due to the singular structure of this new topology.

	Let $X$ be a normed space, $0<p<+\infty$ and $(f_k)_k$ a sequence in $X^\ast$. The sequence $(f_k)_k$ is said to be weak$^\ast$ $\textrm{w}_p-$ summable to $f\in X^\ast$ if for every $x\in X$, $f_k(x)\overset{\textrm{w}_p}{\to}f(x)$, that is,
$$
	\lim_n\frac{1}{n}\displaystyle\sum_{k=1}^n|f_k(x) - f(x)|^p=0, 
$$
and we will write $(f_k)\overset{w^\ast-\textrm{w}_p}{\to}f$ and $f=w^\ast-\textrm{w}_p-\lim_n f_n$.

Let $\sum_i f_i$ be a series in the dual space $X^\ast$ of a Banach space $X$, $0<p<+\infty$. We now consider the space of $w^\ast-\textrm{w}_p$ summability defined by:
$$
	S_{w^\ast-\textrm{w}_p}\left(\sum_i f_i\right)=\left\{(a_i)_i\in\ell_{\infty}: \sum_{i} a_i f_i\mbox{ is }w^\ast-\textrm{w}_p\mbox{ summable}\right\}
$$
endowed with the supremum norm. 

\begin{theorem}\label{Sw*-wpcompleto}
	Let $X$ be a normed space and $\sum f_i$ a series in $X^\ast$. Let us consider the following statements:
    \medskip
    \begin{enumerate}[(1)]
    	\item $\sum_i f_i$ is a weakly unconditionally Cauchy series (wuc).  
    	\item $S_{w^\ast-\textrm{w}_p}(\sum_i f_i)=\ell_{\infty}$. 
        \item If $x\in X$ and $M\subset\mathbb{N}$, then $\sum_{i\in M}f_i(x)$ is $\emph{\textrm{w}}_p$ convergent.
    \end{enumerate}
    Then $(1)\Rightarrow(2)\Rightarrow(3)$, and if $X$ is barrelled, then $(3)\Rightarrow(1)$.
\end{theorem}

{\it Proof. } If $(a_i)_i\in\ell_{\infty}$, then the series $\sum_i a_if_i$ is $w^\ast$ convergent in $X^\ast$, that is, there exists $f\in X^{\ast}$ such that $\sum_{i=1}^n a_if_i\overset{w^\ast}{\to}f$. This implies that for every $x\in X$, $\sum_i f_i(x)=f(x)$, and it is easily shown that $\sum_{i=1}^n a_if_i(x)\overset{\textrm{w}_p}{\to}f(x)$, which implies that $(a_i)_i\in S_{w^\ast-\textrm{w}_p}$.  

The implication $(2)\Rightarrow (3)$ follows directly.

Now, if $X$ is barrelled, let us define 
$$ E=\left\{\sum_{i=1}^n a_i f_i: n\in\mathbb{N}, |a_i|\leq 1\right\}.
$$ 
In order to prove $(3)\Rightarrow(1)$, it is sufficient to show that $E$ is pointwise bounded. Suppose on the contrary that there exists $x_0\in X$ such that $\sum_i|f_i(x_0)|$ diverges. If $M^+=\{i\in\mathbb{N}:f_i(x_0)\geq0\}$ and $M^-=\{i\in\mathbb{N}:f_i(x_0)<0\}$, then either $\sum_{i\in M^+}f_i(x_0)$ diverges or $\sum_{i\in M^-}(-f_i)(x_0)$ diverges. Then, by applying Proposition \ref{propoStolz} and Proposition \ref{propo}, we obtain that the series is not $\emph{\textrm{w}}_p$ convergent, which is a contradiction with $(3)$.$\hfill\Box$



\begin{thebibliography}{10}

\bibitem{b3}
Abdullah Alotaibi and M.~Mursaleen.
\newblock Korovkin type approximation theorems via lacunary equistatistical convergence.
\newblock {\em Filomat}, 30(13):3641--3647, 2016.

\bibitem{boos}
Johann Boos.
\newblock {\em Classical and modern methods in summability}.
\newblock Oxford Mathematical Monographs. Oxford University Press, Oxford, 2000.
\newblock Assisted by Peter Cass, Oxford Science Publications.

\bibitem{connorganichevkadets}
J.~Connor, M.~Ganichev, and V.~Kadets.
\newblock A characterization of {B}anach spaces with separable duals via weak statistical convergence.
\newblock {\em J. Math. Anal. Appl.}, 244(1):251--261, 2000.

\bibitem{b4}
J.~Connor and E.~Sava\c~s.
\newblock Lacunary statistical and sliding window convergence for measurable functions.
\newblock {\em Acta Math. Hungar.}, 145(2):416--432, 2015.

\bibitem{Connor88}
J.~S. Connor.
\newblock The statistical and strong {$p$}-{C}es\`aro convergence of sequences.
\newblock {\em Analysis}, 8(1-2):47--63, 1988.

\bibitem{fekete}
M.~Fekete.
\newblock Viszg{\'a}latok a fourier-sorokr{\'o}l (research on fourier series).
\newblock {\em Math. {\'e}s term{\'e}sz}, 34:759--786, 1916.

\bibitem{Fridy}
J.~A. Fridy.
\newblock On statistical convergence.
\newblock {\em Analysis}, 5(4):301--313, 1985.

\bibitem{b5}
J.~A. Fridy and C.~Orhan.
\newblock Lacunary statistical convergence.
\newblock {\em Pacific J. Math.}, 160(1):43--51, 1993.

\bibitem{hardy}
G.~H. Hardy and J.~E. Littlewood.
\newblock Sur la s{\'e}rie de fourier d'une fonction {\'a} carr{\'e} sommable.
\newblock {\em CR Acad. Sci. Paris}, 156:1307--1309, 1913.

\bibitem{kolk}
Enno Kolk.
\newblock The statistical convergence in {B}anach spaces.
\newblock {\em Tartu \"Ul. Toimetised}, (928):41--52, 1991.

\bibitem{b2}
M.~Mursaleen, S.~Debnath, and D.~Rakshit.
\newblock {$I$}-statistical limit superior and {$I$}-statistical limit inferior.
\newblock {\em Filomat}, 31(7):2103--2108, 2017.

\bibitem{b1}
Mohammad Mursaleen and Cemal Belen.
\newblock On statistical lacunary summability of double sequences.
\newblock {\em Filomat}, 28(2):231--239, 2014.

\bibitem{Aizpuru2000}
F.~J. P\'erez-Fern\'andez, F.~Ben\'\i~tez Trujillo, and A.~Aizpuru.
\newblock Characterizations of completeness of normed spaces through weakly unconditionally {C}auchy series.
\newblock {\em Czechoslovak Math. J.}, 50(125)(4):889--896, 2000.

\bibitem{b7}
Ekrem Sava\c~s.
\newblock Generalized asymptotically {$I$}-lacunary statistical equivalent of order {$\alpha$} for sequences of sets.
\newblock {\em Filomat}, 31(6):1507--1514, 2017.

\bibitem{b6}
Ekrem Sava\c~s.
\newblock {$\mathcal{I}_\lambda$}-statistically convergent functions of order {$\alpha$}.
\newblock {\em Filomat}, 31(2):523--528, 2017.

\bibitem{b8}
Ekrem Sava\c~s.
\newblock On {$\mathcal{I}$}-lacunary statistical convergence of weight {$g$} of sequences of sets.
\newblock {\em Filomat}, 31(16):5315--5322, 2017.

\bibitem{zeller-beekmann}
K.~Zeller and W.~Beekmann.
\newblock {\em Theorie der {L}imitierungsverfahren}.
\newblock Zweite, erweiterte und verbesserte Auflage. Ergebnisse der Mathematik und ihrer Grenzgebiete, Band 15. Springer-Verlag, Berlin-New York, 1970.

\bibitem{zygmund}
Antoni Zygmund.
\newblock {\em Trigonometrical series}.
\newblock Dover Publications, New York, 1955.

\end{thebibliography}

\end{document}